\documentclass[12pt, a4paper, reqno]{amsart}
\usepackage{amsfonts,amsmath,latexsym,verbatim,amscd,mathrsfs,color,array}
\usepackage[colorlinks=true]{hyperref}
\usepackage{amssymb,amsthm,graphicx}
\usepackage{geometry}
\usepackage{bbm}
\usepackage{pdfsync}
\usepackage{epstopdf}
\usepackage{cite}
\usepackage{multirow}
\usepackage[font=bf,aboveskip=15pt]{caption}
\usepackage[toc,page]{appendix}
\usepackage{hyperref}
\hypersetup{
    colorlinks=true,
    linkcolor=red,
    citecolor=blue,
    urlcolor=black
}
\allowdisplaybreaks[3]

\usepackage{titlesec}
\titleformat{\section}
  {\large\bfseries\raggedright} 
  {\thesection.} 
  {0.4em} 
  {} 

\usepackage{fancyhdr}

\newtheorem{theo}{Theorem}[section]
\newtheorem{lemm}{Lemma}[section]

\newtheorem{defi}{Definition}[section]

\numberwithin{equation}{section}

\usepackage{marvosym} 

\usepackage{xcolor} 

\begin{document}

\title{\large Complete Classification of the Symmetry Group of $L_p$-Minkowski Problem on the Sphere}
\author{Huan-Jie Chen$^1$ and Shi-Zhong Du$^1$}

\thanks{\Letter\enspace Shi-Zhong Du (szdu@stu.edu.cn) \enspace\&\enspace Huan-Jie Chen (19hjchen@stu.edu.cn)}
\thanks{$^1$\enspace\enspace The Department of Mathematics, Shantou University, Shantou 515063, P. R. China.}
\thanks{*\enspace\enspace This work is partially supported by Natural Science Foundation of China (12171299)}

\maketitle

\renewenvironment{abstract}
 {\small\noindent\textbf{Abstract}\enspace}
 {\par\vspace{0.8em}}

\begin{abstract}
   In Convex Geometry, a core topic is the $L_p$-Minkowski problem
   \begin{equation}\label{e0.1}
     \det(\nabla^2h+hI)=fh^{p-1}, \ \ \forall X\in{\mathbb{S}}^n, \ \ \forall p\in \mathbb{R}
   \end{equation}
  of Monge-Amp\`{e}re type. By the transformation $u(x)=h(X)\sqrt{1+|x|^2}$ and semi-spherical projection, equation \eqref{e0.1} can be reformulated by the Monge-Amp\`{e}re type equation
   \begin{equation}\label{e0.2}
     \det D^2u=(1+|x|^2)^{-\frac{p+n+1}{2}}u^{p-1}, \ \ \forall x\in{\mathbb{R}}^n, \ \ \forall p\in \mathbb{R}
   \end{equation}
  on the Euclidean space. In this paper, we will firstly determine the symmetric groups of $n$-dimensional fully nonlinear equation \eqref{e0.2} without asymptotic growth assumption. After proving several key resolution lemmas, we thus completely classify the symmetric groups of the $L_p$-Minkowski problem. Our method develops the Lie theory to fully nonlinear PDEs in Convex Geometry.
\end{abstract}


\vspace{8pt}


\maketitle


\vspace{15pt}

\section{Introduction}

\noindent The Minkowski problem is to determine a convex body with prescribed Gauss curvature or other similar geometric data. It plays a central role in the theory of convex bodies. Various Minkowski problems \cite{Al,CL1,CY,D,GLL,GM,HLYZ} have been studied especially after Lutwak \cite{L2,L3}, who proposes two variants of the Brunn-Minkowski theory including the dual Brunn-Minkowski theory and the  $L_p$ Brunn-Minkowski theory. Furthermore, there are singular cases such as the logarithmic Minkowski problem \cite{BLYZ2,CLZ,Z1} and the centro-affine Minkowski problem \cite{BLYZ1,CW,L1}.

Given a convex body $K\subset{\mathbb{R}}^{n+1}$ containing the origin, for each $X\in{\mathbb{S}}^n$, let $r(X)$ be a point on the boundary $\partial K$ whose unit outer normal vector is $X$. The support function $h:{\mathbb{S}}^n\rightarrow \mathbb{R}$ of $K$ is defined to be
   \begin{equation}
     h(X)\equiv{r}(X)\cdot X, \ \ \forall X\in{\mathbb{S}}^n.
   \end{equation}
For a uniformly convex $C^2$ body, the matrix $[\nabla^2_{ij}h+h\delta_{ij}]$ is positive definite, where $[\nabla^2_{ij}h]$ stands for the Hessian matrix of $h(X)$ acting on an orthonormal frame $\{e_i\}_{i=1}^n$ of ${\mathbb{S}}^n$ and $\delta_{ij}$ denotes the Kronecker delta symbol. Conversely, any $C^2$ function $h(X)$ satisfying $[\nabla^2_{ij}h+h\delta_{ij}]>0$ determines a uniformly convex $C^2$ body $K_h$. Direct calculation shows the standard surface measure of $\partial K$ is given by
   $$
    dS\equiv\frac{1}{n+1}\det(\nabla^2_{ij}h+h\delta_{ij})d {\mathcal{H}}^{n}|_{{\mathbb{S}}^n}
   $$
for $n$-dimensional Hausdorff measure $d {\mathcal{H}}^{n}$. It is well known that classical Minkowski problem looks for a convex body such that its standard surface measure matches a given Radon measure on ${\mathbb{S}}^n$. In \cite{L3}, Lutwak introduces the $L_p$-surface measure $dS_p\equiv h^{1-p}dS$ on $\partial K$. The corresponding $L_p$-Minkowski problem is to look for a convex body whose $L_p$-surface measure is equal to a prescribed function.

Parallel to the classical Minkowski problem, the $L_p$-Minkowski problem amounts to the solvability of fully nonlinear equation
   \begin{equation}\label{e1.2}
     \det(\nabla^2_{ij}h+h\delta_{ij})=fh^{p-1}, \ \ \forall X\in{\mathbb{S}}^n, \ \ \forall p\in \mathbb{R}
   \end{equation}
in the smooth category. This problem has been widely studied in the past, see for examples \cite{BBCY,C,GLW,LW1,LW2,HLW,Z2} and especially the corner stone paper Chou-Wang \cite{CW}. Regarding the existence of solutions to the $L_p$-Minkowski problem, we have summarize the following conclusion based on known results.

\begin{theo}
Considering the $L_p$-Minkowski problem \eqref{e1.2}, one has

 (1) when $p>n+1$, there exists a unique positive solution in $C^{2,\alpha}({\mathbb{S}}^n)$ for each positive function $f\in C^\alpha({\mathbb{S}}^n), \alpha\in(0,1)$;

 (2) when $p=n+1$, there exists a unique pair $(h,\lambda)$ for $0<h\in C^{2+\alpha}({\mathbb{S}}^n)$ and $\lambda\in{\mathbb{R}}_+$ satisfying
   \begin{equation*}
     \det(\nabla^2_{ij}h+h\delta_{ij})=\lambda fh^n
   \end{equation*}
 for each positive function $f\in C^\alpha({\mathbb{S}}^n), \alpha\in(0,1)$;

 (3) when $1<p<n+1$, \eqref{e1.2} has a generalized non-negative solution in the sense of Aleksandrov for each $f\in L^\infty({\mathbb{S}}^n)$ with $f\geq c_0$, where $c_0$ is some positive constant;

 (4) when $p=1$, \eqref{e1.2} reduces to the classical Minkowski problem and it has a solution if and only if the positive function $f$ satisfies the Kazdan-Warner condition
 $$\int_{{\mathbb{S}}^n}Xfd\sigma=0;$$

 (5) when $-n-1<p<1$, \eqref{e1.2} has a generalized non-negative solution in the sense of Aleksandrov for each $f\in L^\infty({\mathbb{S}}^n)$ with $f\geq c_0$, where $c_0$ is some positive constant. In addition, if $-n-1<p\leq-n+1$ and $f\in C^\alpha({\mathbb{S}}^n)$ for some $\alpha\in(0,1)$, this special solution is positive and in $C^{2,\alpha}({\mathbb{S}}^n)$;

 (6) when $p=-n-1$, \eqref{e1.2} admits infinitely many not affine-equivalent $C^{2,\alpha}$ solutions for each positive function $f\in C^\infty({\mathbb{S}}^n)$;

 (7) when $p<-n-1$, \eqref{e1.2} admits a uniformly convex, smooth and positive solution for each positive function $f\in C^\infty({\mathbb{S}}^n)$.
\end{theo}

When the function $f\equiv1$, equation \eqref{e1.2} is reduced to
   \begin{equation}\label{e1.3}
     \det(\nabla^2_{ij}h+h\delta_{ij})=h^{p-1}, \ \ \forall X\in{\mathbb{S}}^n, \ \ \forall p\in \mathbb{R}.
   \end{equation}
The equation \eqref{e1.3} on the semi-sphere can be projected into the Monge-Amp\`{e}re type equation
   \begin{equation}\label{e1.4}
     \det D^2u=(1+|x|^2)^{-\frac{p+n+1}{2}}u^{p-1}, \ \ \forall x\in{\mathbb{R}}^n, \ \ \forall p\in \mathbb{R}
   \end{equation}
in Euclidean space by $u(x)=h(X)\sqrt{1+|x|^2}$ and the semi-spherical projection which is defined in Section \ref{sec3}.

On the other orientation, a classical method to find symmetry reductions of partial differential equations is the Lie group method. This method has been extensively investigated and widely applied to differential equations in Euclidean space as well as to geometric equations in Riemannian manifold. For example, Wo-Yang-Wang \cite{WYW} obtained group invariant solutions for centro-affine invariant flow by classifying symmetric groups and optimal systems within the framework of Lie's theory. The readers may also refer to other papers \cite{CL2,CDG,CQ,CM,F,OS} and the classical book \cite{O} for the related topic. In the frame of symmetry group, Lie's theory illustrates all transformations preserving the equations invariant, and then classifies the corresponding group invariant solutions. Although there are vast literatures on the $L_p$-Minkowski problem, there are few papers discuss the symmetric groups of this problem using the Lie's theory. This is the main purpose of this paper to classify the symmetry group of \eqref{e1.3} completely, using some algebraic lemmas and resolution lemmas.

By acting the prolongation vector field $pr^{(2)}\overrightarrow{v}$ of the infinitesimal generator $\overrightarrow{v}$ to the $n$-dimensional fully nonlinear projective equation \eqref{e1.4}, and subsequently solving the resulting complex system of equations for the coefficients of $\overrightarrow{v}$ obtained through term-by-term comparison, we establish the following classification result.
\begin{theo}\label{t1.2}
  The symmetry group of the equation \eqref{e1.4} are generated by
    \begin{eqnarray*}
     g_1(\varepsilon) &: & (x,u)\to(y,v)=(R_\varepsilon x, u),\ R_\varepsilon\in SO(n), \\
     g_2(\varepsilon) &: & (x,u)\to(y,v)=(x, \varepsilon u), \\
     g_3^i(\varepsilon) &: & (x,u)\to(y,v)=\frac{(x^1,\cdots,\sin\varepsilon+x^i\cos\varepsilon,\cdots,x^n,u)}{\cos \varepsilon-x^i\sin \varepsilon}
    \end{eqnarray*}
 for $p=n+1$, and $g_1(\varepsilon)$, $g_3^i(\varepsilon)$,
    \begin{eqnarray*}
     g_4(\varepsilon) &: & (x,u)\to(y,v)=(x, u+\varepsilon), \\
     g_5^i(\varepsilon) &: & (x,u)\to(y,v)=(x, u+\varepsilon x^i)
    \end{eqnarray*}
 for $p=1$, and
    \begin{eqnarray*}
     g_6(\varepsilon) &: & (x,u)\to(y,v)=(A_\varepsilon x, u),\ A_\varepsilon\in SL(n), \\
     g^i_7(\varepsilon) &: & (x,u)\to(y,v)=(x^1, \cdots, e^{\varepsilon}x^i, \cdots, x^n, e^{\frac{\varepsilon}{n+1}}u), \\
     g^i_8(\varepsilon) &: & (x,u)\to(y,v)=(x^1, \cdots, x^i+\varepsilon, \cdots, x^n, u), \\
     g^i_9(\varepsilon) &: & (x,u)\to(y,v)=\frac{(x,u)}{1-\varepsilon x^i}
    \end{eqnarray*}
 for $p=-n-1$, and $g_1(\varepsilon)$, $g_3^i(\varepsilon)$ for $p\neq n+1,1,-n-1$, where $i=1,2,\cdots,n$.
\end{theo}

Since there is no asymptotic growth assumption is imposed, our classification result for the equation \eqref{e1.4} is complete. Using Theorem \ref{t1.2} and several resolution lemmas, we thus derive the following characterization result regarding the symmetry group of the $L_p$-Minkowski problem.
\begin{theo}\label{t1.3}
 Considering the symmetry group of the equation \eqref{e1.3}, one has

 (1) when $p=n+1$, it is corresponding to the rotation and scaling transformations of convex bodies;

 (2) when $p=1$, it is corresponding to the rotation and translation transformations of convex bodies;

 (3) when $p=-n-1$, it is corresponding to the centro-affine transformation of convex bodies;

 (4) when $p\neq n+1,1,-n-1$, it is corresponding to the rotation transformation of convex bodies.
\end{theo}

The contents of the paper are organized as follows. Firstly, we will recall some basical facts about the Lie's theory especially the key prolongation formula for PDEs in Section \ref{sec2}. Secondly, we classify the symmetry group of equation \eqref{e1.4} and thus complete the proof of Theorem \ref{t1.2} in Section \ref{sec3}. Thirdly, with the help of the resolution lemmas shown in Sections \ref{sec4}-\ref{sec6} for the cases $p=n+1,1,-n-1$, we complete the proof of the main result Theorem \ref{t1.3} in Section \ref{sec7} by recovering the transformations of the solutions for the projective equation to the transformations on convex bodies.

\vspace{20pt}

\section{Preliminary to Lie's theory on partial differential equations}\label{sec2}

\noindent Let us recall some facts of Lie's theory to partial differential equation. We suggest the reader refer to the classical book \cite{O} by Peter Olver for a detailed introduction on the symmetry group of PDEs. At the beginning, we present the concept of a symmetry group of the equation.

\begin{defi}
Let $F: {\mathbb{R}}^n\times{\mathbb{R}}^m\times{\mathbb{R}}^{mn}\times{\mathbb{R}}^{mn^2}\times\cdots\times{\mathbb{R}}^{mn^k} \to{\mathbb{R}}^l$ be a mapping from ${\mathbb{R}}^{n+m\frac{1-n^{k+1}}{1-n}}$ to ${\mathbb{R}}^l$. A symmetry group of the equation
$$F(x, u, Du, \cdots, D^ku)=0$$
is a local group of transformations $G$ defined on an open subset
$$M\subset\big\{(x,u)\in\mathbb{R}^n\times\mathbb{R}^m\big\}$$
of the space of independent and dependent variables for the equation with the property that if $u=f(x)$ is a solution to the equation $F=0$, then whenever $g\circ f$ is defined for $g\in G$, we have $u=g\circ f(x)$ is also a solution.
\end{defi}

If one considering a partial differential equation
   \begin{equation}\label{e2.1}
     F(x, u, Du, D^2u)=0
   \end{equation}
of second order, and suppose that
   $$
    \overrightarrow{v}=\xi^i(x,u)\frac{\partial}{\partial x^i}+\phi(x,u)\frac{\partial}{\partial u}
   $$
is an infinitesimal generator of one-parameter group action $g(\varepsilon)\ (\varepsilon\in\mathbb{R})$ of equation \eqref{e2.1}, we have the second prolongation of $\overrightarrow{v}$ is the vector field
   \begin{equation}\label{e2.2}
    pr^{(2)}\overrightarrow{v}=\xi^i\frac{\partial}{\partial x^i}+\phi\frac{\partial}{\partial u}+\phi^i\frac{\partial}{\partial u_i}+\phi^{ij}\frac{\partial}{\partial u_{ij}}
   \end{equation}
by prolongation formula in \cite{O} (Theorem $2.36$, Page $110$), where
    \begin{eqnarray*}
       \phi^i&=&D_i(\phi-\xi^ju_j)+\xi^ju_{ij}\\
       &=&\phi_i+\phi_uu_i-(\xi^j_i+\xi^j_uu_i)u_j,\\
      \phi^{ij}&=&D_{ij}(\phi-\xi^au_a)+\xi^au_{ija}\\
       &=&D_j\big\{\phi_i+\phi_uu_i-(\xi^a_i+\xi^a_uu_i)u_a-\xi^au_{ia}\big\}+\xi^au_{ija}\\
       &=&\phi_{ij}+\phi_{iu}u_j+(\phi_{uj}+\phi_{uu}u_j)u_i+\phi_uu_{ij}\\
       &&-\big[(\xi^a_{ij}+\xi^a_{iu}u_j)+(\xi^a_{uj}+\xi^a_{uu}u_j)u_i+\xi^a_uu_{ij}\big]u_a\\
       &&-(\xi^a_i+\xi^a_uu_i)u_{ja}-(\xi^a_j+\xi^a_uu_j)u_{ia}.
    \end{eqnarray*}
Moreover, $g(\cdot)$ is an one-parameter symmetry group of equation \eqref{e2.1} if and only if
   \begin{equation}
      pr^{(2)}\overrightarrow{v}F(x,u,Du,D^2u)=0
   \end{equation}
holds for any $u^{(2)}\equiv(u, Du, D^2u)$ satisfying
   \begin{equation}
     F(x, u^{(2)})=0,
   \end{equation}
where $x, u, Du, D^2u$ are regarded as independent variables as usually.

\vspace{15pt}

\section{Projecting of the $L_p$-Minkowski problem from ${\mathbb{S}}^n_-$ to ${\mathbb{R}}^n$}\label{sec3}

\noindent Letting $X=(X',X_{n+1})\in{\mathbb{R}}^n\times{\mathbb{R}}$ be the coordinate representation of the southern semi-sphere ${\mathbb{S}}^n_-$, we have the semi-spherical projection $T:{\mathbb{S}}^n_-\rightarrow \mathbb{R}^n$ is given by
$$T(X',X_{n+1})\equiv (x,-1)=\bigg(-\frac{X'}{X_{n+1}},-1\bigg), \  \forall (X',X_{n+1})\in {\mathbb{S}}^n_-.$$
Regarding $x=(x^1,\cdots,x^n)\in \mathbb{R}^n$ as local coordinates of ${\mathbb{S}}^n_-$, one can express $X$ in terms of $x$ by
    \begin{equation}
     X=\bigg(\frac{x}{\sqrt{1+|x|^2}},-\frac{1}{\sqrt{1+|x|^2}}\bigg).
    \end{equation}
Upon this local coordinates, the induced Riemannian metric $g$ of ${\mathbb{S}}^n_-$ is
   \begin{eqnarray*}
    g_{ij}&\equiv &\left(\frac{\partial X}{\partial x^i},\frac{\partial X}{\partial x^j}\right)\\
    &=&\frac{1}{1+|x|^2}\bigg(\delta_{ij}-\frac{x^ix^j}{1+|x|^2}\bigg),
   \end{eqnarray*}
where $i,j=1,2,\cdots,n$. Moreover, the inverse metric and determinate of $g_{ij}$ are given by
  $$
  \begin{cases}
   g^{ij}=(1+|x|^2)(\delta_{ij}+x^ix^j),\\
   \det g_{ij}=(1+|x|^2)^{-(n+1)}
  \end{cases}
  $$
respectively. Thus, the Christoffel symbol equals
   \begin{eqnarray*}
    \Gamma^k_{ij}&\equiv&\frac{1}{2}g^{kl}\bigg(\frac{\partial g_{jl}}{\partial x^i}+\frac{\partial g_{il}}{\partial x^j}-\frac{\partial g_{ij}}{\partial x^l}\bigg)\\
    &=&-\frac{\delta_{jk}x^i+\delta_{ik}x^j}{1+|x|^2}
   \end{eqnarray*}
for different combinatory indices $i, j, k$. Under an orthogonal frame $\{e_i\}_{i=1}^n$ of the sphere, the $L_p$-Minkowski problem boils down to solve the fully nonlinear equation
   \begin{equation}\label{e3.2}
   \det (\nabla^2_{ij}h+h\delta_{ij})=fh^{p-1},
   \end{equation}
where $h(X)$ is the support function and $f$ is the density of a given Borel measure $d\mu=fd{\mathcal{H}}^n|_{{\mathbb{S}}^n}$ for $n$ dimensional Hausdorff measure $d{\mathcal{H}}^n$. In the local coordinates defined as above, \eqref{e3.2} can also be written as
   \begin{equation}\label{e3.3}
   \frac{\det(\nabla^2_{ij}h+hg_{ij})}{\det g_{ij}}=fh^{p-1}.
   \end{equation}
Assume that the support function $h(X(x))$ and its corresponding projective function $u(x)$ satisfy the relationship
   \begin{equation}
   h(X(x))=\frac{u(x)}{\sqrt{1+|x|^2}},
   \end{equation}
then direct calculation shows that the Hessian matrix of $h$ is given by
   \begin{eqnarray*}
   \nabla^2_{ij}h&\equiv&\frac{\partial^2h}{\partial x^i\partial x^j}-\Gamma^k_{ij}\frac{\partial h}{\partial x^k}\\
   &=&\frac{u_{ij}}{\sqrt{1+|x|^2}}-\frac{u\delta_{ij}}{(1+|x|^2)^\frac{3}{2}}+ \frac{ux^ix^j}{(1+|x|^2)^{\frac{5}{2}}}.
   \end{eqnarray*}
Summing as above, if one considers the case $f\equiv1$, \eqref{e3.3} will be reduced to the Monge-Amp\`{e}re type equation
   \begin{equation}\label{e3.5}
   \det D^2u=(1+|x|^2)^{-\frac{p+n+1}{2}}u^{p-1},\ \ x\in \mathbb{R}^n.
   \end{equation}
At this point, we have projected the $L_p$-Minkowski problem on the semi-sphere into the equation \eqref{e3.5} on the Euclidean space.

Now, let us prove the classification result Theorem \ref{t1.2} for the symmetry group of the projective equation \eqref{e3.5} by applying the Lie's theory. The steps of Lie's theory is firstly to find the conditions on infinitesimal generator $\overrightarrow{v}$ of an one parameter symmetry group of the equation. This is a linear condition on $\overrightarrow{v}$, or prolongation of $\overrightarrow{v}$. After comparing the coefficients of like terms on both sides of the equation, we obtain all infinitesimal generators $\overrightarrow{v}$. Then, using the one-to-one corresponding of Lie algebra with Lie group, one can recover all symmetric groups of the equation.

\vspace{4pt}

\noindent\textbf{Proof of Theorem \ref{t1.2}:}
Setting the infinitesimal generator of the one parameter symmetry group by
  $$
   \overrightarrow{v}=\xi^i(x,u)\frac{\partial}{\partial x^i}+\phi(x,u)\frac{\partial}{\partial u},
  $$
we have the prolongation formula for $\overrightarrow{v}$ up to second order is given by
   $$
     pr^{(2)}\overrightarrow{v}=\xi^i\frac{\partial}{\partial x^i}+\phi\frac{\partial}{\partial u}+\phi^{i}\frac{\partial}{\partial u_{i}}+\phi^{ij}\frac{\partial}{\partial u_{ij}},
   $$
where the coefficients $\phi^{i}$ and $\phi^{ij}$ have been presented in Section \ref{sec2}. Define
$$\Phi(x,u^{(2)}):=\det D^2u-(1+|x|^2)^{-\frac{p+n+1}{2}}u^{p-1},$$
then the derivatives of $\Phi(x,u^{(2)})$ are
\begin{equation*}
\left\{
\begin{aligned}
&\mathcal{X}^i:=\frac{\partial}{\partial x^i}\big[\Phi(x,u^{(2)})\big]=(p+n+1)(1+|x|^2)^{-\frac{p+n+3}{2}}u^{p-1}x^i,\\
&\mathcal{U}:=\frac{\partial}{\partial u}\big[\Phi(x,u^{(2)})\big]=(1-p)(1+|x|^2)^{-\frac{p+n+1}{2}}u^{p-2},\\
&\mathcal{D}_{ij}:=\frac{\partial}{\partial u_{ij}}\big[\Phi(x,u^{(2)})\big]=U^{ij},
\end{aligned}
\right.
\end{equation*}
where $U^{ij}$ denotes the co-factor matrix of $u_{ij}$. Therefore, $\overrightarrow{v}$ is an infinitesimal generator of a one-parameter symmetry group of \eqref{e3.5} if and only if
\begin{equation}\label{e3.6}
  \xi^i\mathcal{X}^i+\phi\thinspace\mathcal{U}+\phi^{ij}\mathcal{D}_{ij}=0.
\end{equation}
Notice that there holds $u_{ik}U^{kj}=\delta_{ij}\det D^2u$ for any $i,j$, then
   $$
     \phi^{ij}\mathcal{D}_{ij}=\eta^{ij}U^{ij}+\Big[n\phi_u -2\textstyle\sum\limits_k\xi^k_k-(n+2)\xi^k_uu_k\Big](1+|x|^2)^{- \frac{p+n+1}{2}}u^{p-1},
   $$
where
$$\eta^{ij}=\phi_{ij}+\phi_{iu}u_j+(\phi_{uj}+\phi_{uu}u_j)u_i- \big[\xi^k_{ij}+\xi^k_{iu}u_j+(\xi^k_{uj}+\xi^k_{uu}u_j)u_i\big]u_k.$$
Checking the first order derivative term, we immediately obtain that $\xi^k_uu_k=0$, or equivalent to $\xi^i_u=0$ for all $i$. For another, by comparing the second order derivative terms on both sides of \eqref{e3.6}, one obtains that $\phi_{uu}=0$, $\phi_{ij}=0$ and
$$\phi_{iu}u_j+\phi_{uj}u_i-\xi^k_{ij}u_k=0,\ \forall i,j.$$
Suppose that $\xi^i$ and $\phi$ take the forms
\begin{equation}
\begin{cases}
\xi^i(x,u)=a_kx^kx^i+A^i_kx^k+B^i,\\
\phi(x,u)=(a_kx^k+b)u+c_kx^k+d,
\end{cases}
\end{equation}
where $a_i,b,c_i,d,A^i_j,B^i$ all are constants. Then
$$\textstyle\sum\limits_k\xi^k_k=(n+1)\textstyle\sum\limits_ka_kx^k+\textstyle\sum\limits_kA^k_k,$$
and \eqref{e3.6} can be simplified by
\begin{equation}\label{e3.8}
  (p+n+1)\mathcal{T}_1-(p-1)\mathcal{T}_2-\mathcal{T}_3=0,
\end{equation}
where the terms $\mathcal{T}_1,\mathcal{T}_2$ and $\mathcal{T}_3$ are given by
\begin{eqnarray*}
\mathcal{T}_1&=&(a_kx^k|x|^2+A^i_kx^kx^i+B^ix^i)u,\\
\mathcal{T}_2&=&(1+|x|^2)\big[(a_kx^k+b)u+c_kx^k+d\big],\\
\mathcal{T}_3&=&(1+|x|^2)\Big[(n+2)a_kx^k+2\textstyle\sum\limits_kA^k_k-nb\Big]u.
\end{eqnarray*}
Next, our situations are divided into four cases:

\noindent\textbf{Case 1: $p=-n-1$.} \eqref{e3.8} changes to
   \begin{eqnarray*}
     \Big[2(n+1)b-2\textstyle\sum\limits_kA^k_k\Big]u+(n+2)(c_kx^k+d)=0,
   \end{eqnarray*}
then $(n+1)b=\sum_kA^k_k$ and $c_i=d=0$ for any $i$. Hence, one concludes that
\begin{equation}
\begin{cases}
\xi^i(x,u)=a_kx^kx^i+A^i_kx^k+B^i,\\
\phi(x,u)=(a_kx^k+b)u.
\end{cases}
\end{equation}
And the Lie algebra of infinitesimal generators are spanned by the vector fields
  \begin{eqnarray*}
    \overrightarrow{v}^i_1&=&\partial_{x^i},\\
    \overrightarrow{v}_2&=&\textstyle\sum\limits_{i\neq j}A^i_jx^j\partial_{x^i},\\
    \overrightarrow{v}_3&=&\textstyle\sum\limits_iA^i_ix^i\partial_{x^i}+bu\partial_u,\\
    \overrightarrow{v}^i_4&=&x^ix^j\partial_{x^j}+x^iu\partial_u.
  \end{eqnarray*}

\noindent\textbf{Case 2: $p=1$.} By combining the like terms, \eqref{e3.8} can be simplified by
   $$
    (n+2)A^i_kx^kx^i+\Big(nb-2\textstyle\sum\limits_kA^k_k\Big)|x|^2 +(n+2)(B^ix^i-a_kx^k)+nb-2\textstyle\sum\limits_kA^k_k=0,
   $$
then $b=0$, $A^i_j=-A^j_i$ and $B^i=a_i$ for all $i,j$. Hence, it follows that
\begin{equation}
\begin{cases}
\xi^i(x,u)=a_kx^kx^i+A^i_kx^k+a_i,\\
\phi(x,u)=a_kx^ku+c_kx^k+d.
\end{cases}
\end{equation}
And the Lie algebra of infinitesimal generators are spanned by the vector fields
  \begin{eqnarray*}
    \overrightarrow{v}_1&=&\partial_u,\\
    \overrightarrow{v}^i_2&=&x^i\partial_u,\\
    \overrightarrow{v}_3&=&A^i_jx^j\partial_{x^i},\\
    \overrightarrow{v}^i_4&=&x^ix^j\partial_{x^j}+\partial_{x^i}+x^iu\partial_u.
  \end{eqnarray*}

When $p\neq -n-1,1$, by comparing the terms that do not contain $u$ on both sides of \eqref{e3.8}, one derives that $c_i=d=0$ for any $i$. Then \eqref{e3.8} equivalent to
\begin{equation}\label{3.11}
(p+n+1)\mathcal{T}_4-(p-n-1)\mathcal{T}_5-\mathcal{T}_6=0,
\end{equation}
where the terms $\mathcal{T}_4,\mathcal{T}_5$ and $\mathcal{T}_6$ are given by
\begin{eqnarray*}
\mathcal{T}_4&=&A^i_kx^kx^i+B^ix^i-a_kx^k,\\
\mathcal{T}_5&=&b(1+|x|^2),\\
\mathcal{T}_6&=&2\textstyle\sum\limits_kA^k_k(1+|x|^2).
\end{eqnarray*}

\noindent\textbf{Case 3: $p=n+1$.} \eqref{3.11} changes to
   \begin{eqnarray*}
     (n+1)(A^i_kx^kx^i+B^ix^i-a_kx^k)-\textstyle\sum\limits_kA^k_k(1+|x|^2)=0,
   \end{eqnarray*}
it yields that
\begin{equation}
\begin{cases}
\xi^i(x,u)=a_kx^kx^i+A^i_kx^k+a_i,\\
\phi(x,u)=(a_kx^k+b)u,
\end{cases}
\end{equation}
where $A^i_j=-A^j_i$ for all $i,j$. And the infinitesimal generators are spanned by
  \begin{eqnarray*}
    \overrightarrow{v}_1&=&u\partial_u,\\
    \overrightarrow{v}_2&=&A^i_jx^j\partial_{x^i},\\
    \overrightarrow{v}^i_3&=&x^ix^j\partial_{x^j}+\partial_{x^i}+x^iu\partial_u.
  \end{eqnarray*}

\noindent\textbf{Case 4: $p\neq -n-1,1,n+1$.} Similarly, according to \eqref{3.11}, we have
\begin{equation}
\begin{cases}
\xi^i(x,u)=a_kx^kx^i+A^i_kx^k+a_i,\\
\phi(x,u)=a_kx^ku,
\end{cases}
\end{equation}
where $A^i_j=-A^j_i$ for all $i,j$. And the infinitesimal generators are spanned by
  \begin{eqnarray*}
    \overrightarrow{v}_1&=&A^i_jx^j\partial_{x^i},\\
    \overrightarrow{v}^i_2&=&x^ix^j\partial_{x^j}+\partial_{x^i}+x^iu\partial_u.
  \end{eqnarray*}
Thus, we have completed the proof of Theorem \ref{t1.2}. \hfill $\Box$

\vspace{15pt}

\section{Resolution transformations on convex body for $p=n+1$}\label{sec4}

\noindent As in Section \ref{sec3}, we project the equation \eqref{e1.3} on ${\mathbb{S}}^n_-$ to obtain an equation \eqref{e1.4} on ${\mathbb{R}}^n$. After applying Lie's theory to \eqref{e1.4}, we derive all transformations that preserving \eqref{e1.4} invariant. To completely classify the symmetry group of the $L_p$-Minkowski problem, we need to recover symmetry group of the projective equation \eqref{e1.4} on convex bodies.

To begin with, let us set up some notations for simplicity. For any group action
$$g(\varepsilon):(x,u)\rightarrow (y,v),\ \varepsilon\in\mathbb{R},\ x,y\in \mathbb{R}^n,$$
we denote the original convex body as $K_1$ and the convex body corresponding to the transformation under the group action as $K_2$. Assume that their unit outer normal vectors are represented as
\begin{equation}
\left\{
\begin{aligned}
X(x)&=\big(X_1(x),\cdots,X_{n+1}(x)\big),\\
Y(y)&=\big(Y_1(y),\cdots,Y_{n+1}(y)\big)\\
\end{aligned}
\right.
\end{equation}
respectively, or equivalent to
\begin{equation}\label{e.normal1}
\left\{
\begin{aligned}
X(x)&=\bigg(\frac{x}{\sqrt{1+|x|^2}},-\frac{1}{\sqrt{1+|x|^2}}\bigg),\\
Y(y)&=\bigg(\frac{y}{\sqrt{1+|y|^2}},-\frac{1}{\sqrt{1+|y|^2}}\bigg).
\end{aligned}
\right.
\end{equation}
Additionally, we employ the notations
\begin{equation}
\left\{
\begin{aligned}
Z(X)&=\big(Z_1(X),\cdots,Z_{n+1}(X)\big),\\
W(Y)&=\big(W_1(Y),\cdots,W_{n+1}(Y)\big)
\end{aligned}
\right.
\end{equation}
to denote the position vectors of $K_1$ and $K_2$ respectively. We claim that all vectors are column vectors, and we omit the transposition symbol for convenience. Furthermore, the support functions of $K_1$ and $K_2$ are given by
\begin{equation}\label{e4.4}
\left\{
\begin{aligned}
h(X(x))&=\frac{u(x)}{\sqrt{1+|x|^2}},\\
q(Y(y))&=\frac{v(y)}{\sqrt{1+|y|^2}}.
\end{aligned}
\right.
\end{equation}

By the definition of support function, we immediately obtain the following lemma, which illustrates the relationship between the rotation transformation of a convex body and the transformation in its support function.
\begin{lemm}\label{lem4.1}
For two convex bodies $K_1$ and $K_2$, suppose that their unit outer normal vectors are related by $Y=SX$ for some matrix $S\in SO(n+1)$. Then, the relation $W=SZ$ holds if and only if their support functions satisfy $q(Y)=h(X)$.
\end{lemm}
Firstly, we apply Lemma \ref{lem4.1} to the group action
\begin{equation}
g_1(\varepsilon)\ :\ (x,u)\rightarrow (y,v)=(R_\varepsilon x,u),\ R_\varepsilon\in SO(n).
\end{equation}
In this case, the unit outer normal vector of the transformed convex body $K_2$ is
   \begin{eqnarray}\nonumber
   Y&=&\bigg(\frac{R_\varepsilon x}{\sqrt{1+|x|^2}},-\frac{1}{\sqrt{1+|x|^2}}\bigg)\\
   &=&S_\varepsilon X,
   \end{eqnarray}
where the matrix $S_\varepsilon$ is given by
     \begin{equation*}
       S_\varepsilon=\left[
       \begin{array}{cc}
         R_\varepsilon & 0 \\
         0 & 1
       \end{array}
       \right]\in SO(n+1).
     \end{equation*}
By \eqref{e4.4}, it is derived that the support function of $K_2$ equals
   \begin{eqnarray}\nonumber
   q(Y(y))&=&\frac{u(R_\varepsilon x)}{\sqrt{1+|x|^2}}\\
   &=&h(X(R_\varepsilon x)).
   \end{eqnarray}
Thus, according to Lemma \ref{lem4.1}, this group action locally corresponds to a rotation transformation of convex bodies that preserve the $(n+1)$-th coordinate.

The following lemma elucidates the relationship between the scaling transformation of a convex body and the transformation in its support function.
\begin{lemm}\label{lem4.2}
For two convex bodies $K_1$ and $K_2$, the relation $W_i=k_iZ_i,i=1,2,\cdots,n+1$, where $k_i\in\mathbb{R}$, holds if and only if there exist unit outer normal vectors $X$ and $Y$, such that their support functions satisfy
   \begin{equation}\label{e4.8}
   q(Y)=h(X)\sqrt{k_1^2Y_1^2+\cdots+k_{n+1}^2Y_{n+1}^2}.
   \end{equation}
In particular, the relation $W=kZ,k\in\mathbb{R}$ holds if and only if $Y=X$ and
   \begin{equation}
   q(Y)=kh(X).
   \end{equation}
\end{lemm}
\noindent\textbf{Proof:} First prove the necessity. Let $\varphi(Z)=0$ be the representation function of convex body $K_1$, then the unit outer normal vectors of $K_1$ and $K_2$ are represented as
   \begin{equation}
   \left\{
   \begin{aligned}
   X&=\frac{(\varphi_{Z_1}, \cdots,\varphi_{Z_{n+1}})}{\sqrt{\varphi_{Z_1}^2 +\cdots+\varphi_{Z_{n+1}}^2}},\\
   Y&=\frac{(k_1^{-1}\varphi_{Z_1},\cdots,k_{n+1}^{-1}\varphi_{Z_{n+1}})}{\sqrt{ k_1^{-2}\varphi_{Z_1}^2 +\cdots+k_{n+1}^{-2}\varphi_{Z_{n+1}}^2}}
   \end{aligned}
   \right.
   \end{equation}
respectively. Therefore, the support functions of $K_1$ and $K_2$ are given by
   \begin{eqnarray*}
   h(X)&=&\frac{Z_1\varphi_{Z_1}+\cdots+Z_{n+1}\varphi_{Z_{n+1}}}{\sqrt{\varphi_{Z_1}^2 +\cdots+\varphi_{Z_{n+1}}^2}},\\
   q(Y)&=&\frac{Z_1\varphi_{Z_1}+\cdots+Z_{n+1}\varphi_{Z_{n+1}}}{\sqrt{k_1^{-2}\varphi_{Z_1}^2+ \cdots+k_{n+1}^{-2}\varphi_{Z_{n+1}}^2}}\\
   &=&h(X)\sqrt{\frac{\varphi_{Z_1}^2+\cdots+\varphi_{Z_{n+1}}^2}{k_1^{-2}\varphi_{Z_1}^2+ \cdots+k_{n+1}^{-2}\varphi_{Z_{n+1}}^2}}
   \end{eqnarray*}
respectively, which implies \eqref{e4.8}. For another, since the convex body and its support function uniquely determine each other, the converse statement also holds true. \hfill $\Box$

\vspace{5pt}

Secondly, we apply Lemma \ref{lem4.2} to the group action
\begin{equation}
g_2(\varepsilon)\ :\ (x,u)\rightarrow (y,v)=(x,\varepsilon u).
\end{equation}
At this time, the unit outer normal vectors of the convex bodies $K_1$ and $K_2$ satisfy $X=Y$. And the support function of $K_2$ is given by
   \begin{eqnarray}\nonumber
   q(Y(y))&=&\frac{\varepsilon u(x)}{\sqrt{1+|x|^2}}\\
   &=&\varepsilon h(X(x)).
   \end{eqnarray}
By Lemma \ref{lem4.2}, we conclude that this group action corresponds to the scaling transformation of the entire convex bodies by a factor of $k=\varepsilon$.

Thirdly, we consider the group action
\begin{equation}
g_3^i(\varepsilon)\ :\ (x,u)\rightarrow (y,v)=\frac{(x^1,\cdots,\sin\varepsilon+x^i\cos\varepsilon,\cdots,x^n,u)}{\cos \varepsilon-x^i\sin \varepsilon},\ \forall i.
\end{equation}
Direct calculation shows that the unit outer normal vector of $K_2$ is
   \begin{eqnarray}\nonumber
   Y&=&\frac{(x^1,\cdots,x^i\cos\varepsilon+\sin\varepsilon,\cdots,x^n,x^i\sin\varepsilon-\cos\varepsilon)}{\sqrt{1+|x|^2}}\\
   &=&O_\varepsilon X,
   \end{eqnarray}
where the matrix $O_\varepsilon$ is given by
     \begin{equation}\label{e4.15}
       O_\varepsilon=\left[
       \begin{array}{cccc}
         I_{i-1} & 0 & 0 & 0\\
         0 & \cos\varepsilon & 0 & -\sin\varepsilon\\
         0 & 0 & I_{n-i} & 0\\
         0 & \sin\varepsilon & 0 & \cos\varepsilon
       \end{array}
       \right]\in SO(n+1).
     \end{equation}
Based on the relationship between $x$ and $y$, we have
\begin{equation*}
\left\{
\begin{aligned}
   x^i&=\frac{y^i\cos \varepsilon-\sin\varepsilon}{y^i\sin\varepsilon+\cos\varepsilon}, \\
   x^j&=\frac{y^j}{y^i\sin\varepsilon+\cos\varepsilon}, \ \ \forall j\neq i.
\end{aligned}
\right.
\end{equation*}
By the semi-spherical projection, one can express $y$ in terms of $Y$ by
\begin{equation}\label{e4.16}
     y=-\frac{(Y_1,\cdots,Y_n)}{Y_{n+1}},
\end{equation}
and it yields that
   $$
    |y|^2=\frac{1}{Y_{n+1}^2}-1.
   $$
Therefore, the support function of the transformed convex body $K_2$ is given by
   \begin{eqnarray}\nonumber
     q_3(Y)&=&\frac{h(X)}{\cos\varepsilon-x^i\sin\varepsilon}\sqrt{\frac{1+|x|^2}{1+|y|^2}}\\
     &=&h(X).
   \end{eqnarray}
By Lemma \ref{lem4.1}, we conclude that this group action locally corresponds to a rotation transformation \eqref{e4.15} of convex bodies.

In summary, we have successfully recovered the group actions $g_1(\varepsilon)$, $g_2(\varepsilon)$ and $g_3^i(\varepsilon)$ on convex bodies. Consequently, we derive the following resolution lemma.
\begin{lemm}\label{lem4.3}
 For the symmetry group of equation \eqref{e1.4}, when $p=n+1$, all group actions locally correspond to certain compositions of rotation and scaling transformations of convex bodies.
\end{lemm}

\vspace{15pt}

\section{Resolution transformations on convex body for $p=1$}\label{sec5}

\noindent To recover the group actions $g_4(\varepsilon)$ and $g^i_5(\varepsilon)$ on convex body, we need the following lemma, which describes the relationship between the translation of a convex body and the transformation in its support function.
\begin{lemm}\label{lem5.1}
For two convex bodies $K_1$ and $K_2$, the relation $W_i=Z_i+b_i,b_i\in\mathbb{R}$, where $i=1,2,\cdots,n+1$, holds if and only if $X=Y$ and
   \begin{equation}\label{e5.1}
   q(Y)=h(X)+b_1Y_1+\cdots+b_{n+1}Y_{n+1}.
   \end{equation}
\end{lemm}
\noindent\textbf{Proof:} Using the notations from the proof of Lemma \ref{lem4.2}, we have
   $$Y=\frac{(\varphi_{Z_1}, \cdots,\varphi_{Z_{n+1}})}{\sqrt{\varphi_{Z_1}^2+\cdots+\varphi_{Z_{n+1}}^2}}$$
and thus $X=Y$. The support function of $K_2$ is given by
   \begin{equation*}
   q(Y)=\frac{(Z_1+b_1)\varphi_{Z_1}+\cdots+(Z_{n+1}+b_{n+1})\varphi_{Z_{n+1}}}{\sqrt{\varphi_{Z_1}^2 +\cdots+\varphi_{Z_{n+1}}^2}},
   \end{equation*}
which implies \eqref{e5.1} and the necessity has been proven. On the other hand, the sufficiency also holds true. \hfill $\Box$

\vspace{5pt}

Firstly, we apply Lemma \ref{lem5.1} to the group action
\begin{equation}
g_4(\varepsilon)\ :\ (x,u)\rightarrow (y,v)=(x,u+\varepsilon).
\end{equation}
For this case, the unit outer normal vectors of the convex bodies $K_1$ and $K_2$ satisfy $X=Y$. And using \eqref{e4.16}, we obtain that the support function of $K_2$ is given by
   \begin{eqnarray}\nonumber
   q(Y)&=&h(X)+\frac{\varepsilon}{\sqrt{1+|y|^2}}\\
   &=&h(X)+\varepsilon Y_{n+1}.
   \end{eqnarray}
By Lemma \ref{lem5.1}, this group action locally corresponds to a translation transformation of the convex bodies along the $(n+1)$-th axis by $\varepsilon$ units.

Secondly, we consider the group action
\begin{equation}
g_5^i(\varepsilon)\ :\ (x,u)\rightarrow (y,v)=(x,u+\varepsilon x^i),\ \forall i.
\end{equation}
At this time, the unit outer normal vectors of $K_1$ and $K_2$ satisfy $X=Y$. By \eqref{e4.16}, one derives that the support function of $K_2$ equals
   \begin{eqnarray}\nonumber
   q_5(Y)&=&h(X)+\frac{\varepsilon y^i}{\sqrt{1+|y|^2}}\\
   &=&h(X)-\varepsilon Y_i.
   \end{eqnarray}
According to Lemma \ref{lem5.1}, we conclude that this group action locally corresponds to a translation transformation of the convex bodies alone the negative direction of the $i$-th axis by $\varepsilon$ units.

Combined with the analysis of the group actions $g_1(\varepsilon)$ and $g_3^i(\varepsilon)$ in Section \ref{sec4}, this yields the following resolution lemma for the case $p=1$.
\begin{lemm}\label{lem5.2}
 For the symmetry group of equation \eqref{e1.4}, when $p=1$, all group actions locally correspond to certain compositions of rotation and translation transformations of convex bodies.
\end{lemm}

\vspace{15pt}

\section{Resolution transformations on convex body for $p=-n-1$}\label{sec6}

\noindent Let us start with a basic algebraic decomposition lemma, it tell us that the element of special linear group $SL(n)$ can be decomposed using compositions of special orthogonal matrices and diagonal matrices.

\begin{lemm}\label{lem6.1}
Let $A\in SL(n)$, then $A$ admits a matrix decomposition of the form
   \begin{equation}
    A=P\mathrm{diag}(\lambda_1,\cdots,\lambda_n)Q,
   \end{equation}
where $P,Q\in SO(n)$ and $\lambda_1,\cdots,\lambda_n>0$ satisfies $\lambda_1\cdots\lambda_n=1$.
\end{lemm}
\noindent\textbf{Proof:} Note that $A^TA$ is a positive definite matrix since $A$ is invertible. So, there exist an orthogonal matrix $B$ and positive constants $\mu_1,\cdots,\mu_n$ with $\mu_1\cdots\mu_n=1$, such that
$$A^TA=B^T\mathrm{diag}(\mu_1,\cdots,\mu_n)B.$$
Let $C=B^T\mathrm{diag}(\sqrt{\mu_1},\cdots,\sqrt{\mu_n})B$, then $C^2=A^TA$ and $C$ is a positive definite matrix.
Since $C^2=C^TC$, we have
$$(AC^{-1})^T(AC^{-1})=I.$$
Thus, $D=AC^{-1}$ is an orthogonal matrix and
$$A=DC=DQ^T\mathrm{diag}(\lambda_1,\cdots,\lambda_n)Q,$$
where $Q\in SO(n)$. Let $P=DQ^T$, we obtain the desired decomposition.  \hfill $\Box$

\vspace{5pt}

Firstly, we consider the group action
\begin{equation}
g_6(\varepsilon)\ :\ (x,u)\rightarrow (y,v)=(A_\varepsilon x,u),\ A_\varepsilon\in SL(n).
\end{equation}
By Lemma \ref{lem6.1}, the elements in $SL(n)$ can be decomposed as the product of some orthogonal matrices and diagonal matrices with a determinant equal to $1$. Therefore, we can divide $g_6(\varepsilon)$ into two parts $g_1(\varepsilon)$ and
\begin{equation}
\widetilde{g}_6^{ij}(\mu)\ :\ (x,u)\rightarrow (y,v)=(x^1,\cdots,\mu^{-1}x^i,\cdots,\mu x^j,\cdots,x^n,u),\ \forall i,j,
\end{equation}
where $\mu>0$. For the group action $\widetilde{g}_6^{ij}$, it is obvious that the unit outer normal vector of $K_2$ is given by
   \begin{eqnarray}\nonumber
   Y&=&\frac{(x^1,\cdots,\mu^{-1}x^i,\cdots,\mu x^j,\cdots,x^n,-1)}{\sqrt{(x^1)^2+\cdots+\mu^{-2}(x^i)^2+\cdots+\mu^2(x^j)^2 +\cdots+(x^n)^2+1}}\\
   &=&\frac{(X_1,\cdots,\mu^{-1}X_i,\cdots,\mu X_j,\cdots,X_n,X_{n+1})}{\sqrt{X_1^2+\cdots+\mu^{-2}X_i^2+\cdots+\mu^2X_j^2 +\cdots+X_n^2+X_{n+1}^2}}.
   \end{eqnarray}
Using \eqref{e4.16}, we get that the support function of the transformed convex body is
   \begin{eqnarray}\nonumber
   \widetilde{q}_6(Y)&=&h(X)\sqrt{\frac{1+|x|^2}{1+|y|^2}}\\ \nonumber
   &=&h(X)\sqrt{1+(\mu^{-2}-1)Y_i^2+(\mu^2-1)Y_j^2}\\
   &=&h(X)\sqrt{\mu^{-2}Y_i^2+\mu^2Y_j^2+{\textstyle\sum}_{s\neq i,j}Y_s^2}.
   \end{eqnarray}
According to Lemma \ref{lem4.2}, one obtains that the scaling factors are
$$k_i=\mu^{-1},~k_j=\mu,\ k_s=1,\ \forall s\neq i,j.$$
Thus, this group action locally corresponds to the scaling transformation of convex bodies and the scaling factors satisfy $k_1\cdots k_{n+1}=1$.

Secondly, we consider the group action
\begin{equation}
g_7^i(\mu)\ :\ (x,u)\rightarrow (y,v)=(x^1, \cdots, \mu x^i, \cdots, x^n, \mu^\frac{1}{n+1}u),\ \forall i,
\end{equation}
where we set $\mu=e^\varepsilon>0$ for simplicity. At this time, the unit outer normal vector of $K_2$ is given by
   \begin{eqnarray}\nonumber
   Y&=&\frac{(x^1,\cdots,\mu x^i,\cdots,x^n,-1)}{\sqrt{(x^1)^2+\cdots+\mu^2(x^i)^2+ \cdots+(x^n)^2+1}}\\ \nonumber
   &=&\frac{\big(\mu^{-\frac{1}{n+1}}x^1,\cdots,\mu^{\frac{n}{n+1}} x^i,\cdots,\mu^{-\frac{1}{n+1}}x^n,-\mu^{-\frac{1}{n+1}}\big)}{\sqrt{ \mu^{-\frac{2}{n+1}}(x^1)^2+\cdots+\mu^{\frac{2n}{n+1}}(x^i)^2 +\cdots+\mu^{-\frac{2}{n+1}}(x^n)^2+\mu^{-\frac{2}{n+1}}}}\\
   &=&\frac{\big(\mu^{-\frac{1}{n+1}}X_1,\cdots,\mu^{\frac{n}{n+1}} X_i,\cdots,\mu^{-\frac{1}{n+1}}X_n,\mu^{-\frac{1}{n+1}}X_{n+1}\big)}{\sqrt{ \mu^{-\frac{2}{n+1}}X_1^2+\cdots+\mu^{\frac{2n}{n+1}}X_i^2 +\cdots+\mu^{-\frac{2}{n+1}}X_n^2+\mu^{-\frac{2}{n+1}}X_{n+1}^2}}.
   \end{eqnarray}
For another, by \eqref{e4.16}, we derive that the support function of $K_2$ equals
   \begin{eqnarray}\nonumber
   q_7(Y)&=&\mu^\frac{1}{n+1}h(X)\sqrt{\frac{1+|x|^2}{1+|y|^2}} \\ \nonumber
   &=&\mu^\frac{1}{n+1}h(X) \sqrt{1+(\mu^{-2}-1)Y_i^2}\\
   &=&\mu^\frac{1}{n+1}h(X) \sqrt{Y_1^2+\cdots+\mu^{-2}Y_i^2+\cdots+Y^2_{n+1}}.
   \end{eqnarray}
By Lemma \ref{lem4.2}, one obtains that the scaling factors are determined as
$$k_i=\mu^{-\frac{n}{n+1}},\ k_j=\mu^\frac{1}{n+1},\ \forall j\neq i.$$
Hence, we conclude that this group action locally corresponds to the scaling transformation of convex bodies and the scaling factors satisfy $k_1\cdots k_{n+1}=1$.

To recover the group action $g^i_8(\varepsilon)$ on a convex body, we require the following lemma, which establishes the relationship between the certain centro-affine transformation of a convex body and the corresponding transformation of its support function.
\begin{lemm}\label{lem6.2}
For two convex bodies $K_1$ and $K_2$, the relation $W=H_\varepsilon Z$,
     \begin{equation}\label{e6.9}
       H_\varepsilon=\left[
       \begin{array}{cccc}
         I_{i-1} & 0 & 0 & 0\\
         0 & 1 & 0 & 0\\
         0 & 0 & I_{n-i} & 0\\
         0 & \varepsilon & 0 & 1
       \end{array}
       \right]\in SL(n+1)
     \end{equation}
holds if and only if there exist the unit outer normal vectors $X$ and $Y$, such that their support functions satisfy
   \begin{equation}\label{e6.10}
   q(Y)=h(X)\sqrt{1+2\varepsilon Y_iY_{n+1}+\varepsilon^2Y_{n+1}^2}.
   \end{equation}
\end{lemm}
\noindent\textbf{Proof:} Once again, use the notations from the proof of Lemma \ref{lem4.2}. According to the relationship between $W$ and $Z$, one can express $Z$ in terms of $W$ by
\begin{equation*}
\left\{
\begin{aligned}
   &Z_j=W_j, \ j=1,2,\cdots,n,\\
   &Z_{n+1}=W_{n+1}-\varepsilon W_i.
\end{aligned}
\right.
\end{equation*}
Then $Y$ can be represented as
\begin{equation}\label{e6.11}
   Y=\frac{\big(\varphi_{Z_1},\cdots,\varphi_{Z_i}-\varepsilon\varphi_{Z_{n+1}},\cdots, \varphi_{Z_{n+1}}\big)}{\sqrt{\varphi_{Z_1}^2 +\cdots+(\varphi_{Z_i}-\varepsilon \varphi_{Z_{n+1}})^2+\cdots+\varphi_{Z_{n+1}}^2}}.
\end{equation}
And the support function of transformed convex body $K_2$ is given by
   \begin{eqnarray}\nonumber\label{e6.12}
   q(Y)&=&\frac{Z_1\varphi_{Z_1}+\cdots+Z_i(\varphi_{Z_i}-\varepsilon\varphi_{Z_{n+1}}) +\cdots+Z_n\varphi_{Z_n}+\varphi_{Z_{n+1}}(\varepsilon Z_i+Z_{n+1})}{\sqrt{\varphi_{Z_1}^2 +\cdots+(\varphi_{Z_i}-\varepsilon \varphi_{Z_{n+1}})^2+\cdots+\varphi_{Z_{n+1}}^2}}\\
   &=&h(X)\sqrt{\frac{\varphi_{Z_1}^2+\cdots+\varphi_{Z_{n+1}}^2}{\varphi_{Z_1}^2 +\cdots+(\varphi_{Z_i}-\varepsilon \varphi_{Z_{n+1}})^2+\cdots+\varphi_{Z_{n+1}}^2}}.
   \end{eqnarray}
From \eqref{e6.11}, we obtain
   \begin{equation*}
   Y_i^2=\frac{\varphi_{Z_i}^2-2\varepsilon\varphi_{Z_i}\varphi_{Z_{n+1}}+\varepsilon^2\varphi_{Z_{n+1}}}{\varphi_{Z_1}^2+\cdots+(\varphi_{Z_i}-\varepsilon \varphi_{Z_{n+1}})^2+\cdots+\varphi_{Z_{n+1}}^2},
   \end{equation*}
and solving this equation yields that
\begin{equation}\label{e6.13}
\frac{\varphi_{Z_i}^2}{\varphi_{Z_1}^2+\cdots+(\varphi_{Z_i}-\varepsilon \varphi_{Z_{n+1}})^2+\cdots+\varphi_{Z_{n+1}}^2}=Y_i^2+2\varepsilon Y_iY_{n+1}+\varepsilon^2Y_{n+1}^2.
\end{equation}
Substituting \eqref{e6.13} into \eqref{e6.12}, we derive
   \begin{equation*}
   q(Y)=h(X)\sqrt{{\textstyle\sum}_{k\neq i}Y_k^2+Y_i^2+2\varepsilon Y_iY_{n+1}+\varepsilon^2Y_{n+1}^2},
   \end{equation*}
which implies the identity \eqref{e6.10}. Besides, the sufficiency also valid. \hfill $\Box$

\vspace{5pt}

Thirdly, we apply Lemma \ref{lem6.2} to the group action
\begin{equation}
g^i_8(\varepsilon)\ :\ (x,u)\rightarrow (y,v)=(x^1, \cdots, x^i+\varepsilon, \cdots, x^n, u),\ \forall i.
\end{equation}
For this case, the unit outer normal vector of $K_2$ is
   \begin{eqnarray}\nonumber
   Y&=&\frac{(x^1,\cdots, x^i+\varepsilon,\cdots,x^n,-1)}{\sqrt{(x^1)^2+\cdots+(x^i+\varepsilon)^2+ \cdots+(x^n)^2+1}}\\
   &=&\frac{(X_1,\cdots,X_i-\varepsilon X_{n+1},\cdots,X_n,X_{n+1})}{\sqrt{1-2\varepsilon X_iX_{n+1}+\varepsilon^2X_{n+1}^2}}.
   \end{eqnarray}
From \eqref{e4.16}, the support function of $K_2$ is given by
   \begin{eqnarray}\nonumber\label{e6.14}
   q_8(Y)&=&h(X)\sqrt{\frac{1+|x|^2}{1+|y|^2}}\\
   &=&h(X)\sqrt{1+2\varepsilon Y_iY_{n+1}+\varepsilon^2Y_{n+1}^2}.
   \end{eqnarray}
Applying Lemma \ref{lem6.2}, we conclude that this group action locally corresponds to the centro-affine transformation \eqref{e6.9} of convex bodies.

Similarly, to recover the group action $g^i_9(\varepsilon)$ on a convex body, we require the following result, which is analogous to Lemma \ref{lem6.2} and we omit the proof for brevity.
\begin{lemm}\label{lem6.3}
For two convex bodies $K_1$ and $K_2$, the relation $W=Q_\varepsilon Z$,
     \begin{equation}\label{e6.17}
       Q_\varepsilon=\left[
       \begin{array}{cccc}
         I_{i-1} & 0 & 0 & 0\\
         0 & 1 & 0 & -\varepsilon\\
         0 & 0 & I_{n-i} & 0\\
         0 & 0 & 0 & 1
       \end{array}
       \right]\in SL(n+1)
     \end{equation}
holds if and only if there exist the unit outer normal vectors $X$ and $Y$, such that their support functions satisfy
   \begin{equation}
   q(Y)=h(X)\sqrt{1-2\varepsilon Y_iY_{n+1}+\varepsilon^2Y_{n+1}^2}.
   \end{equation}
\end{lemm}
Fourthly, we apply Lemma \ref{lem6.3} to the group action
\begin{equation}
g^i_9(\varepsilon)\ :\ (x,u)\rightarrow (y,v)=\frac{(x,u)}{1-\varepsilon x^i},\ \forall i.
\end{equation}
At this time, the unit outer normal vector of $K_2$ is given by
   \begin{eqnarray}\nonumber
   Y&=&\frac{\big(x,-(1-\varepsilon x^i)\big)}{\sqrt{|x|^2+(1-\varepsilon x^i)^2}}\\
   &=&\frac{(X_1,\cdots,X_n,\varepsilon X_i+X_{n+1})}{\sqrt{1+2\varepsilon X_iX_{n+1}+\varepsilon^2X_i^2}}.
   \end{eqnarray}
Based on the relationship between $x$ and $y$, one can express $x$ in terms of $y$ by
$$x^j=\frac{y^j}{1+\varepsilon y^i},\ j=1,2,\cdots,n.$$
It follows from \eqref{e4.16} that the support function of $K_2$ is
   \begin{eqnarray}\nonumber\label{e6.16}
     q_9(Y)&=&\frac{h(X)}{1-\varepsilon x^i}\sqrt{\frac{1+|x|^2}{1+|y|^2}}\\ \nonumber
     &=&h(X)\sqrt{\frac{(1+\varepsilon y^i)^2+|y|^2}{1+|y|^2}}\\
     &=&h(X)\sqrt{1-2\varepsilon Y_iY_{n+1}+\varepsilon^2Y_i^2}.
   \end{eqnarray}
Applying Lemma \ref{lem6.3}, we conclude that this group action locally also corresponds to the centro-affine transformation \eqref{e6.17} of convex bodies.

Summing as above, we obtain the following resolution lemma for the case $p=-n-1$.
\begin{lemm}\label{lem6.4}
 For the symmetry group of equation \eqref{e1.4}, when $p=-n-1$, all group actions locally correspond to certain centro-affine transformations of convex bodies.
\end{lemm}

\vspace{15pt}

\section{Symmetry group of the $L_p$-Minkowski problem}\label{sec7}

\noindent In Sections \ref{sec3}-\ref{sec6}, we have completely classified all transformations of the projective equation \eqref{e1.4} on ${\mathbb{R}}^n$ and then used resolution lemmas to recover the corresponding transformations on convex bodies. For any $p\in \mathbb{R}$, all transformations that preserve \eqref{e1.4} can be decomposed as compositions of the corresponding group actions. This means that we have classified all transformations of \eqref{e1.3} on the open subset of semi-sphere ${\mathbb{S}}^n_-$. More precisely, we have the following representation lemma.

\begin{lemm}\label{lem7.1}
 Suppose that an open subset $U\subset {\mathbb{S}}^n$ satisfies $diam\thinspace U<\frac{1}{2}diam\thinspace {\mathbb{S}}^n$
 and a mapping $\psi:U\rightarrow{\mathbb{S}}^n$ leaves the equation \eqref{e1.3} invariant, then

 (1) for $p=n+1$, $\psi$ can be expressed as a composition of rotation and scaling transformations;

 (2) for $p=1$, $\psi$ can be expressed as a composition of rotation and translation transformations;

 (3) for $p=-n-1$, $\psi$ can be expressed as a centro-affine transformation;

 (4) for $p\neq n+1,1,-n-1$, $\psi$ can be expressed as a rotation transformation.
\end{lemm}
\noindent\textbf{Proof:} For any $U$, there exists an open subset $V\subset\mathbb{S}^n_-$ and an element $G_0\in SO(n+1)$ such that $U=G_0(V)$. Then $\psi\circ G_0$ maps $V$ to ${\mathbb{S}}^n$ and preserves the geometric equation \eqref{e1.3}. Now, we focus on the case $p=n+1$. By resolution Lemma \ref{lem4.3}, we have
$$\psi\circ G_0=G_1\circ G_2 \circ\cdots\circ G_s,\ s\in \mathbb{Z}_+\ \mathrm{on}\ V,$$
where $G_i\ (i=1,2,\cdots,s)$ are all compositions of rotation and scaling transformations. Notice that rotation and scaling transformations are commutative, so there exist a rotation transformation $R_0$ and a scaling transformation $S$, such that
$$\psi\circ G_0=R_0\circ S.$$
Denote $R=R_0\circ (G_0)^{-1}$, then $R\in SO(n+1)$ and $\psi=R\circ S$.

Similarly, applying resolution Lemma \ref{lem5.2} and noting the commutativity of the rotation and the translation, we establish the case $p=1$. For $p=-n-1$, the conclusion follows directly from the resolution Lemma \ref{lem6.4} and the definition of centro-affine transformation. Finally, the result extends naturally to other values of $p$ by analogous reasoning.  \hfill $\Box$

\vspace{5pt}

By the compactness of ${\mathbb{S}}^n$, there exists a finite open cover $\{\Omega_1,\Omega_2,\cdots,\Omega_k\}$ consisting of sufficiently small subsets, where $k\in\mathbb{N}_+$. According to Lemma \ref{lem7.1}, every local transformation on each $\Omega_i$ can be expressed as a composition of previously established transformations, where $i=1,2,\cdots,k$. Each such local transformation corresponds to a global transformation on ${\mathbb{S}}^n$, and the collection of these global transformations forms the symmetry group of the equation \eqref{e1.3}. For local transformations defined on any two overlapping open subsets of ${\mathbb{S}}^n$, the following compatibility lemma holds.

\begin{lemm}\label{lem7.2}
 For any open subsets $\Omega_1,\Omega_2\subset{\mathbb{S}}^n$ with $\Omega_1\cap\Omega_2\neq \emptyset$, suppose that there exist two local transformations $\psi_1$ on $\Omega_1$ and $\psi_2$ on $\Omega_2$ that leave the equation \eqref{e1.3} invariant. These two local transformations correspond respectively to the global transformations $\rho_1$ and $\rho_2$ on ${\mathbb{S}}^n$. If $\psi_1=\psi_2$ on $\Omega_1\cap\Omega_2$, then $\rho_1\equiv\rho_2$ on ${\mathbb{S}}^n$.
\end{lemm}
\noindent\textbf{Proof:} Without loss of generality, we only consider the case $p=n+1$. Note that
$$\rho_1|_{\Omega_1}\circ\rho_2^{-1}|_{\Omega_2}=Id\ \mathrm{on} \ \Omega_1\cap\Omega_2.$$
By Lemma \ref{lem7.1}, the mapping $\rho_1|_{\Omega_1}\circ\rho_2^{-1}|_{\Omega_2}$ on open subset $\Omega_1\cap\Omega_2\subset{\mathbb{S}}^n$ also corresponding to the composition of rotation and scaling transformations. Hence, $\rho_1\circ\rho_2^{-1}$ must be the identity transformation on ${\mathbb{S}}^n$ and the proof was done.  \hfill $\Box$

\vspace{5pt}

Finally, by combining Lemma \ref{lem7.1} and Lemma \ref{lem7.2}, we directly establish our main result Theorem \ref{t1.3}, thus achieving a complete classification of the symmetry group for the $L_p$-Minkowski problem.

\vspace{20pt}

{\small\noindent\textbf{Acknowledgments}\enspace The author (SZ) would like to express his deepest gratitude to Professors Xi-Ping Zhu, Kai-Seng Chou, Xu-Jia Wang and Neil Trudinger for their constant encouragements and warm-hearted helps. This paper is also dedicated to the memory of Professor Dong-Gao Deng.}

\end{document}